\documentclass{article} 
\usepackage{amssymb} 
\newcommand{\ben}{\begin{enumerate}}
\newcommand{\een}{\end{enumerate}}
\newcommand{\be}{\begin{equation}}
\newcommand{\ee}{\end{equation}}
\newcommand{\bse}{\begin{subequation}}
\newcommand{\ese}{\end{subequation}}
\newcommand{\bea}{\begin{eqnarray}}
\newcommand{\eea}{\end{eqnarray}}
\newcommand{\bc}{\begin{center}}
\newcommand{\ec}{\end{center}}

\newcommand{\ad}{\mbox{ad}\,}
\newcommand{\bb}[1]{{\mathbb #1}}
\newcommand{\cc}[1]{{\cal #1}}

\newtheorem{theorem}{Theorem}

\begin{document}

\author{B.L. Cerchiai,$\strut^{1,2}$ \,
        G. Fiore,$\strut^{3,4}$ \, J. Madore$\strut^{5,2}$ \\\\
        \and
        $\strut^1$Sektion Physik, Ludwig-Maximilian Universit\"at,\\
        Theresienstra\ss e 37, D-80333 M\"unchen
        \and
        $\strut^2$Max-Planck-Institut f\"ur Physik\\
        F\"ohringer Ring 6, D-80805 M\"unchen
        \and
        $\strut^3$Dip. di Matematica e Applicazioni, Fac.  di Ingegneria\\
        Universit\`a di Napoli, V. Claudio 21, 80125 Napoli
        \and
        $\strut^4$I.N.F.N., Sezione di Napoli,\\
        Mostra d'Oltremare, Pad. 19, 80125 Napoli
        \and
        $\strut^5$Laboratoire de Physique Th\'eorique et Hautes Energies\\
        Universit\'e de Paris-Sud, B\^atiment 211, F-91405 Orsay
        }
\date{}
\title{Geometrical techniques for the $N$-dimensional Quantum Euclidean Spaces
${}^*$}

\maketitle
\begin{abstract}
We briefly report our application \cite{CerFioMad00} of a version 
of noncommutative geometry to  the quantum Euclidean space
{\boldmath $R$}${}^N_q$, for any $N\ge 3$;
this space is covariant under the action of the quantum group
$SO_q(N)$, and two covariant differential calculi are known on it.
More precisely, we describe how to construct
in a Cartan `moving-frame formalism' the metric, two covariant derivatives,
the Dirac operator, the frame, the inner derivations dual to the
frame elements, for both of these calculi. 
The components of the frame elements in the basis of differentials
provide a `local realization' of the Faddeev-Reshetikhin-Takhtadjan
generators of $U_q^{\pm}(so(N))$.
\end{abstract}

\vfill{\footnotesize ${}^*$ Talk given by the first author at the 
6th Wigner Symposium, Istanbul, Turkey,\linebreak
16-22 August 1999}
\newpage

\section{Introduction}

The idea that the structure of space-time at short distances may be
well described by a non-commutative geometry
has been appealing since 1947 \cite{Sny47}, because
such noncommutativity might lead to a regularization of the corresponding
field theory (see e.g. \cite{DopFreRob95}, \cite{DicPolSteWes99}, 
\cite{ChoHinMadSte99}). 
Here we apply the formalism of noncommutative geometry~\cite{Con94},
\cite{Wor} to the quantum Euclidean spaces $\bb{R}^N_q$ with 
$N\ge 3$~\cite{FadResTak89}, which are comodule algebras of the quantum 
groups $SO_q(N)$.
To achieve this goal we use a noncommutative generalization~\cite{DimMad96}
of the moving-frame formalism of E.~Cartan. We generalize the results 
which had been previously found~\cite{FioMad99} for $\bb{R}^3_q$.

When $N$ is odd we can follow a scheme similar to the
one developped for $N=3$. For each of the two $SO_q(N)$
covariant differential calculi defined on $\bb{R}^N_q$
we find two torsion free covariant derivatives.
After adding to the algebra a `dilatator' $\Lambda$,
it is possible to construct an (essentially) unique metric,
in such a way that the covariant derivative is compatible with it.
By further enlarging the algebra by the square roots
and inverses of some elements, we are also able to find for each
of the calculi a frame and the derivatives dual to it. 
When $N$ is even, it is necessary
to add also one of the components $K$ of the angular momentum.
Some of the elements we add have vanishing derivative but are none-the-less
noncommutative analogues of non-constant functions. Then their inclusion
can be interpreted as an embedding of the `configuration space' into
part of `phase space'.

In Section~2 we briefly recall the tools of noncommutative
geometry~\cite{Con94} which will be needed. We start with a formal
noncommutative algebra $\cc{A}$ and with a differential calculus
$\Omega^*(\cc{A})$ over it, and define then the concepts of a frame or
`Stehbein'~\cite{DimMad96}, the corresponding metric, covariant
derivative, and generalized Dirac operator~\cite{Con94}.
In Section~3 we shortly review the definiton of
$\bb{R}_q^N$~\cite{FadResTak89} and then the construction
of two $SO_q(N)$-covariant differential
calculi~\cite{CarSchWat91, Ogi92, WesZum90} on $\bb{R}_q^N$
based on the $\hat R$-matrix
formalism. Both yield the de Rham calculus in the commutative limit.
In Section~4 we proceed with the actual construction of the frame over
$\bb{R}^N_q$ and of the inner derivations dual to it. Within this
framework we recover the `Dirac operator',
which had already been found~\cite{Zum97, Ste96}.
We then determine of the metric and the covariant derivatives.

It turns out that the components of the frame 
in the $\xi^i$ basis automatically provide a
`local realization' of $U_q^{\pm}(so(N))$
in the extended algebra of $\bb{R}^N_q$, i.e. 
they satisfy the `RLL' and the `gLL' relations fulfilled by the
$\cc L^{\pm}$~\cite{FadResTak89} generators of
$U_q^{\pm}(so(N))$
and also fulfil the commutation relations of the latter generators
with the coordinates $x^i$. 
In the case of odd $N$ it is possible to 
`glue' them together to get a realization of the whole of
$U_q^{\pm}(so(N))$.

\section{The Cartan moving-frame formalism}

We start by reviewing a noncommutative extension \cite{DimMad96}
of the moving-frame formalism of E. Cartan. The building blocks are a 
noncommutative algebra $\cc{A}$, which in the commutative limit should become
the algebra of functions on a parallelizable manifold $M$,
and a differential calculus $\Omega^*(\cc{A})$ on it, which should reduce
to the ordinary de~Rham differential calculus on $M$ in the same limit. 
The module of the 1-forms $\Omega^1(\cc{A})$ is required to be free of
rank $N$ ($N$ dimension of the manifold), so that it admits a special basis
$\{ \theta^a \}_{1 \le a \le N}$, referred to as `frame' or
`Stehbein', which commutes with the elements of $\cc{A}$:
\be
[f,\theta^a]=0.                                              \label{thetaf}
\ee
We suppose that the basis $\theta^a$ is dual to a set of inner
derivations $e_a=\ad \lambda_a$:
\be
df=e_a f \theta^a=[\lambda_a,f] \theta^a                         \label{df}
\ee
for any $f \in \cc{A}$. Then it is possible to find a formal
`Dirac operator'~\cite{Con94}
\be
\theta = -\lambda_a \theta^a, \mbox{ such that } df = -[\theta,f].\label{dirac}
\ee
We shall require the center $\cc{Z}(\cc{A})$ of $\cc{A}$ to be trivial:
$\cc{Z}(\cc{A})=\bb{C}$, if this condition is not verified, we shall enlarge
the algebra until it does. If we define the (wedge) product $\pi$
in $\Omega^*(\cc{A})$ by relations of the form
\be
\theta^a \theta^b=P^{ab}{}_{cd} \theta^c \otimes \theta^d ,
\qquad P^{ab}{}_{cd} \in \cc{Z}(\cc{A})\label{prod}
\ee 
then the $\lambda_a$ have to satisfy a quadratic relation of the form
\be
2 \lambda_c \lambda_d P^{cd}{}_{ab} - 
\lambda_c F^c{}_{ab} - K_{ab} = 0.   ,
\qquad F^c{}_{ab}, \: K_{ab}  \in \cc{Z}(\cc{A})          \label{manca}
\ee
In the case of the quantum Euclidean spaces $\bb{R}^N_q$
it turns out that $F^c{}_{ab},K_{ab}=0.$

Next, the metric is defined as a nondegenerate $\cc{A}$-bilinear map
\be
g: \Omega^1(\cc{A}) \otimes_{\cc{A}} \Omega^1(\cc{A}) \rightarrow \cc{A}, 
         \label{metric}
\ee
We shall denote
\be
g(\theta^a \otimes \theta^b)=g^{ab}.
\ee
As a further step, a `generalized flip' can be introduced, an
$\cc{A}$-bilinear map
\be
\sigma: \Omega^1(\cc{A}) \otimes_\cc{A} \Omega^1(\cc{A}) \rightarrow
\Omega^1(\cc{A}) \otimes_\cc{A} \Omega^1(\cc{A}),
\qquad \sigma (\theta^a \otimes \theta^b) = 
S^{ab}{}_{cd} \theta^c \otimes \theta^d.                 \label{2.2.4}
\ee
Due to bilinearity
$g^{ab} \in \cc{Z}(\cc{A})=\bb{C}$ and
$S^{ab}{}_{cd}\in \cc{Z}(\cc{A})=\bb{C}$.

The flip is necessary in order to construct a covariant derivative
$D$~\cite{DubMadMasMou96}, i.e. a map
\be
D: \Omega^1(\cc{A}) \rightarrow \Omega^1(\cc{A}) \otimes \Omega^1(\cc{A})
\ee
satisfying a left and right Leibniz rule:
\be
D (f \xi) =  df \otimes \xi + f D\xi, \qquad
D(\xi f) = \sigma (\xi \otimes df) + (D\xi) f. 
\ee
Then the torsion map can be consistently be defined as
\be
\Theta:\Omega^1(\cc{A})\rightarrow\Omega^2(\cc{A}), \qquad \Theta=d-\pi\circ D.
\ee
where bilinearity requires that
\be
\pi \circ (\sigma + 1) = 0.               \label{torsion}
\ee
The curvature map associated to $D$ is defined by
\be
\mbox{Curv} \equiv D^2= \pi_{12} \circ D_2 \circ D, \qquad
\mbox{Curv} (\theta^a) =
- {1 \over 2} R^a{}_{bcd} \theta^c \theta^d \otimes \theta^b.    \label{curv}
\ee
Here $D_2$ is a natural continuation of the map (\ref{2.2.4}) to the 
tensor product $\Omega^1(\cc{A}) \otimes_\cc{A} \Omega^1(\cc{A})$, namely
$D_2(\xi \otimes \eta) = D\xi \otimes \eta + \sigma_{12} (\xi \otimes D\eta)$.

If $F^a{}_{bc}=0$ a torsion-free covariant derivative~\cite{DubMadMasMou96} is
\be
D \xi = - \theta \otimes \xi +\sigma (\xi \otimes \theta).  \label{covdev0}  
\ee
The compatibility of a covariant derivative with the metric
is~\cite{DubMadMasMou95}
\be
g_{23}\circ D_2= d\circ g , \qquad
S^{ae}{}_{df} g^{fg} S^{bc}{}_{eg} = g^{ab} \delta^c_d. \label{met-comp}
\ee
We suppose~\cite{FioMad98} that $\sigma$ satisfies the braid relation.
 
\section{The quantum Euclidean spaces} \label{preliminaries}

In this section some basic results about the $N$-dimensional
quantum Euclidean space $\bb{R}^N_q$ due to~\cite{FadResTak89} are reviewed.
We start with the matrix $\hat{R}$~\cite{FadResTak89} for $SO_q(N,\bb{C})$.
It is a symmetric $N^2 \times N^2$ matrix, and its main property is that it
satisfies the braid relation.
It admits~\cite{FadResTak89} a projector decomposition:
\be
\hat R = q\cc{P}_s - q^{-1}\cc{P}_a + q^{1-N}\cc{P}_t.       \label{projectorR}
\ee
where the $\cc{P}_s$, $\cc{P}_a$, $\cc{P}_t$ are $SO_q(N)$-covariant
$q$-deformations of the symmetric trace-free,
antisymmetric and trace projectors respectively. 
The projector $\cc{P}_t$ projects on a one-dimensional sub-space and 
can be written in the form 
$\cc{P}_t{}_{kl}^{ij} = (g^{sm}g_{sm})^{-1} g^{ij}g_{kl}$. This
leads to the definition of a metric matrix.
It is a $N \times N$ matrix $g_{ij}$, which is a $SO_q(N)$-isotropic
tensor and is a deformation of the ordinary Euclidean metric
$g_{ij}=q^{-\rho_i} \delta_{i,-j}$.
If $n$ is the rank of $SO(N,\bb{C})$, the indices take the values
$i=-n,\ldots,-1,0,1,\ldots n$ for $N$ odd,
and $i=-n,\ldots,-1, 1,\ldots n$ for $N$ even.
Moreover, we have introduced the notation
$\rho_i=(n-\frac{1}{2},\ldots,\frac{1}{2},0,-\frac{1}{2},\ldots,\frac{1}{2}-n)$
for $N$ odd, $(n-1,\ldots,0,0,\ldots,1-n)$ for $N$ even. 

The metric and the braid matrix
satisfy the `$gTT$' relations~\cite{FadResTak89}
\be
g_{il}\,\hat R^{\pm 1}{}^{lh}_{jk} = 
\hat R^{\mp 1}{}^{hl}_{ij}\,g_{lk}, \qquad
g^{il}\,\hat R^{\pm 1}{}_{lh}^{jk} = 
\hat R^{\mp 1}{}_{hl}^{ij}\,g^{lk}.                          \label{gRrel}
\ee

With the help of the projector $\cc{P}_a$,
the $N$-dimensional quantum Euclidean space is defined as the associative
algebra $\bb{R}^N_q$ generated by elements $\{x^i\}_{i=-n,\ldots,n}$ with
relations
\be
\cc{P}_{a}{}^{ij}_{kl} x^k x^l=0.                              \label{xrel}
\ee
or, more explicitly~\cite{Ogi92}
\be
x^i x^j = q x^j x^i \hbox{ for } i<j, i \neq -j, \qquad
[x^i, x^{-i}]=\left\{
\begin{array}{ll}
k \omega_{i-1}^{-1} r^2_{i-1} &\hbox{ for } i > 1\\
0 & \hbox{ for $i=1$, $N$ even,}\\
h r_0^2 & \hbox{ for $i=1$, $N$ odd.}
\end{array}
\right.
\label{explicitx}
\ee
We use the notation here
$\omega_i=q^{\rho_i}+q^{-\rho_i}$, $h=q^{\frac{1}{2}}-q^{-\frac{1}{2}}$,
$k=q^{\frac{1}{2}}-q^{-\frac{1}{2}}$ and
\be
r^2_i=\sum_{k,l=-i}^i g_{kl} x^k x^l, \qquad i \ge 0 \mbox{ for $N$ odd},
\:  i \ge 1 \mbox{ for $N$ even}.
\label{defr}
\ee

For $q \in \bb{R}^+$ a conjugation $(x^i)^*= x^j g_{ji}$
can be defined on $\bb{R}_q^N$ to obtain what is known as quantum real
Euclidean space $\bb{R}^N_q$.

As this will be necessary for the construction of the elements $\lambda_a$,
we enlarge the algebra $\bb{R}_q^N$ with
the real elements $r_i^{\pm 1}=(r_i^2)^{\pm \frac{1}{2}}$, $i=0 \ldots n$.
There is a unique way to postulate their commutation relations with $x^j$
so that the latter give the commutation relations between $r_i^2$ and $x^j$
which  can be drawn from (\ref{xrel}). Namely
\be
x^j r_i=r_i x^j \hbox{ for } |j| \le i, \qquad
x^j r_i=q r_i x^j \hbox{ for } j < -i, \qquad
x^j r_i=q^{-1} r_i x^j \hbox{ for } j > i
\label{xrrel}
\ee
Note that $r\equiv r_n$ turns out to be central.

There are~\cite{CarSchWat91} two quadratic differential calculi 
$\Omega^*(\bb{R}_q^N)$ and $\bar \Omega^*(\bb{R}_q^N)$, which are
covariant with respect to $SO_q(N)$.
\bea
x^i \xi^j = q\,\hat R^{ij}_{kl} \xi^k x^l,\quad &\quad \cc{P}_{s,t}{}_{kl}^{ij}
\xi^k \xi^l = 0
&\mbox{for } \Omega^1(\bb{R}_q^N),               \label{xxirel}        \\
x^i \bar \xi^j = q^{-1}\,\hat R^{-1}{}^{ij}_{kl} \bar\xi^k x^l,    
& \quad \cc{P}_{s,t}{}_{kl}^{ij} \bar \xi^k \bar \xi^l = 0
&\mbox{for } \bar{\Omega}^1(\bb{R}_q^N),           \label{xxistar}
\eea
where $dx^i=\xi^i$ and $\bar{d} x^i=\bar{\xi}^i$.
If a $*$-structure on $\Omega^1(\bb{R}_q^N) \oplus\bar\Omega^1(\bb{R}_q^N)$ 
is defined by setting $(\xi^i)^* = \bar\xi^j g_{ji}$, the two calculi are
seen to be conjugate.

The Dirac operator~\cite{Con94} of \ref{dirac} is easily verified to be
given by~\cite{Ste96,Zum97}.
\bea
\theta&= \omega_n q^{\frac{N}{2}} k^{-1} r^{-2}
g_{ij} x^i \xi^j, &\mbox{for } \Omega^1(\bb{R}_q^N),\\
\bar\theta &= -\omega_n q^{-\frac{N}{2}} k^{-1} r^{-2} g_{ij} x^i \bar\xi^j
&\mbox{for } \Omega^1(\bb{R}_q^N).
\eea
Now, we have the following difficulty.
In Section~2  we required the center of the algebra $\cc{A}$
to be trivial. But the algebra generated by the $x^i$ and $r_j$
has a nontrivial center, therefore the formalism cannot
be directly applied to it. With a general Ansatz of the type
$\theta^a = \theta^a_i \xi^i$, we immediately see that the
condition (\ref{thetaf}) cannot be fulfilled for $r^2$
if $r^2_n\in\cc{Z}(\cc{A})$.
To solve this problem we add to the algebra also a unitary element
$\Lambda$, and its inverse $\Lambda^{-1}$. It is a ``dilatator'', which
satisfies the commutation relations
\be
x^i \Lambda=q \Lambda x^i.                                 \label{xLambda}
\ee
But this is not enough in the case of even $N$.
We have added the elements $r_1^{\pm 1}=(x^{-1} x^1)^{\pm \frac{1}{2}}$
and therefore the center is non trivial even after $\Lambda$
has been added, because the elements $r_1^{-1}x^{\pm 1}$
commute also with $\Lambda$.
We choose to add a Drinfeld-Jimbo generator $K=q^{\frac {H_1}2}$ and its
inverse $K^{-1}$, where $H_1$ belongs to the Cartan subalgebra
of $U_qso(N)$ and represents the component of the angular momentum
in the $(-1,1)$-plane. 
This new element satisfies the commutation relations
\be
K \Lambda=\Lambda k, \qquad K x^{\pm 1}=q^{\pm 1} x^{\pm 1} K,   \qquad
K x^{\pm i}= x^{\pm i} K \hbox{ for } i>1 \label{xkapparel}
\ee

There are many ways to fix the commutation relations of $\Lambda$ with the 
1-forms compatibly with (\ref{xLambda}).  We choose~\cite{Ogi92} 
\be
\xi^i \Lambda=\Lambda \xi^i, \qquad \Lambda d=q d \Lambda. \label{xiLambda}
\ee
This choice has the disadvantage that $\Lambda$ does not satisfy the
Leibniz rule $d(fg)=f dg + (df) g \;\forall f,g \in \bb{R}^N_q$. 
Nevertheless, $\Lambda$ can then be interpreted in a consistent way as an
element of the Heisenberg algebra, because $\Lambda^{-2}$ can be
constructed~\cite{Ogi92} as a simple polynomial in the coordinates and
derivatives, and, moreover, in the next section we shall see that this allows
us to normalize the $\theta^a$ and $\lambda_a$ in such a way as to 
recover $\bb{R}^N$ as geometry in the commutative limit.

As already observed \cite{FioMad99} there are other possibilities, e.g.
we could have set $d \Lambda=0$. This choice, however, is not completely
satisfactory neither, because we would like $df=0$ to hold only for the
analogues of the constant functions and, moreover, with a procedure
similar to the one described previously~\cite{FioMad99} for $N=3$,
we would recover as geometry in the commutative limit $\bb{R} \times
S^{N-1}$ rather than $\bb{R}^N$.

The same discussion which hold for $\Lambda$ can be done to determine the
commutation relations between $K$ and the 1-forms $\xi^i$.
We choose $d K=0$. Then consistency with (\ref{xkapparel}) requires
\be
K \xi^1=q^{\pm 1} \xi^{\pm 1} K, \qquad K \xi^i=\xi^i K \quad
\hbox{ for } i>1.
\ee

\section{Inner derivations, frame, metric and covariant derivatives} 

Now, we would like to proceed with the actual construction of a frame
$\theta^a$ and of the the associated inner derivations $e_a = \ad \lambda_a$
satisfying the conditions in Section~2 for the extended algebra of
$\bb{R}^N_q$. We start with the Ansatz
\be
\theta^a=\theta^a_i \xi^i                                  \label{ansatz}
\ee
for $\theta^a$, but allow the coefficients $\theta^a_i$ to depend
on $\Lambda$. The equation $[r,\theta^a]=0$ fixes the dependence
of the frame on the dilatator to be linear in $\Lambda^{-1}$. 
{From} the duality condition (\ref{df}) one sees immediately that 
the matrices 
\be
e^i_a=[\lambda_a,x^i]                         \label{eia}
\ee
must be inverse to $\theta^a_i$ in the sense that
$e^i_a \theta^a_j=\delta^i_j$, $\theta^a_i e^i_b=\delta^a_b$.

Equation~(\ref{thetaf}) is equivalent to
\be
x^i \theta^a_j=q^{-1} \hat R^{-1}{}^{ki}_{lj} \theta^a_k x^l, \qquad
x^h e_a^i=q \hat R^{hi}_{jk} e_a^j x^k.
\label{xtheta}
\ee
As the $\lambda_a$ have each only one index, while the coefficients
$\theta^a_i$ of the frame have two, it is easier to look for
solutions $\lambda_a$ to the equation
\be
x^h [\lambda_a,x^i]=q \hat{R}^{hi}_{jk} [\lambda_a, x^j] x^k. \label{Rxlambda}
\ee
The inner derivations $e^i_a$ can be easily computed as commutators of
the $\lambda_a$ with the coordinates according to (\ref{eia}), the
matrix $\theta^a_i$ can be recovered as the inverse of $e^i_a$.

The $r$-dependence of $\theta^a$ is fixed by their commutation relations
with $\Lambda$.
We shall require 
\be
[\Lambda, \theta^a]=0 \quad \Longrightarrow \quad [\theta_i^a,\Lambda]=0,
\qquad [e^i_a,\Lambda]=0.
\ee

Our main results are the following theorems~\cite{CerFioMad00}.
\begin{theorem}
\label{theo1}
$N$ independent solutions of Equation~(\ref{Rxlambda}) are given by
\be
\begin{array}{ll}
\lambda_0=\gamma_0 \Lambda (x^0)^{-1} \hbox{ for $N$ odd,} &
\lambda_{\pm 1}=\gamma_{\pm 1} \Lambda (x^{\pm 1})^{-1} K^{\mp 1} 
\hbox{ for $N$ even,} \\[4pt]
\lambda_a=\gamma_a \Lambda r_{|a|}^{-1}r_{|a|-1}^{-1} x^{-a}
&\hbox{ otherwise,} 
\end{array}                                             \label{deflambda}
\ee
where $\gamma_a \in \bb{C}$ are arbitrary normalization constants.
\end{theorem}
This has been proven by a direct computation in~\cite{CerFioMad00}.
The proof is too long to write it here.
\begin{theorem}
\label{theo2}
If the normalization constants in theorem \ref{theo1} satisfy the conditions
\be
\begin{array}{ll}
\gamma_0 = -q^{-\frac{1}{2}} h^{-1} &\quad\mbox{for $N$ odd,} \\[6pt]
\gamma_1 \gamma_{-1}=
\left\{\begin{array}{l}
-q^{-1} h^{-2}\\
k^{-2}
\end{array}\right.
&\quad\!\begin{array}{l}
\mbox{for $N$ odd,} \\
\mbox{for $N$ even,}
\end{array}\\[8pt]
\gamma_a \gamma_{-a} =
-q^{-1} k^{-2} \omega_a \omega_{a-1} &\quad\mbox{for $a>1$}. \nonumber
\end{array}                                               \label{gamma}
\ee
then the elements $\lambda_a$ fulfill among themselves the commutation
relations
\be
\cc{P}_{a}{}^{ab}_{cd} \lambda_a \lambda_b=0  \label{lambdalambda}
\ee
and the matrices $e^i_a=[\lambda_a,x^i]$ satisfy
\bea
RTT-\hbox{relations:}&&  \hat R^{ij}_{kl} e_a^k e_b^l
= e^i_c e^j_d \hat R^{cd}_{ab}  \label{ree} \\
gTT-\hbox{relations:}&&
g^{ab} e^i_a e^j_b=g^{ij} \Lambda^2, \qquad
g_{ij} e^i_a e^j_b=g_{ab} \Lambda^2              \label{gtt}\\
\hbox{normalization:}&&e^0_0 e^0_0=\Lambda^2.              \label{nor}
\eea
\end{theorem}
Again, the proof would be too long and it can been found in \cite{CerFioMad00}.

Let us make some remarks.
It is interesting to note that the commutation relations
(\ref{lambdalambda}) between the $\lambda_a$ are the same as those
(\ref{explicitx}) satisfied by the $x^i$, and therefore the linear and constant
terms in (\ref{manca}) vanish.

The relations (\ref{gamma}) fix only the value of the product
$\gamma_a \gamma_{-a}$. The determination of it can be done e.g.
by applying the $gTT$-relations for $i=-j$.
We see that $\gamma_0^2$ for $N$ odd and $\gamma_1 \gamma_{-1}$ for $N$ even
are positive real numbers, while all the remaining products 
$\gamma_a \gamma_{-a}$ are negative.

Now, an analogous construction can be done for the barred calculus
$\bar\Omega^*(\cc{A})$.
\begin{theorem}
\label{theo3}
$N$ independent solutions of Equation
\be
[\bar \lambda_a,x^i] x^j = 
q^{-1} \hat R^{-1}{}_{kj}^{li} x^l [\bar \lambda_a, x^k].
\ee
are given by
\be
\begin{array}{ll}
\bar\lambda_0=\bar\gamma_0 \Lambda^{-1} (x^0)^{-1} \mbox{ for $N$ odd,} &
\bar\lambda_{\pm 1} = \bar\gamma_{\pm 1} \Lambda^{-1} 
(x^{\pm 1})^{-1} K^{\pm 1} \mbox{ for $N$ even,} \\[6pt]
\bar\lambda_a = 
\bar\gamma_a \Lambda^{-1} r_{|a|}^{-1}r_{|a|-1}^{-1} x^{-a}
&\mbox{otherwise,} 
\end{array}                                         \label{defbarlambda}
\ee
where $\bar \gamma_a \in \bb{C}$ are arbitrary normalization constants.
\end{theorem}

\begin{theorem}
\label{theo4}
If the normalization constants in theorem \ref{theo3} satisfy the conditions
\be
\begin{array}{ll}
\bar \gamma_0 = q^{\frac{1}{2}} h^{-1}
&\quad\mbox{for $N$ odd,} \\[6pt]
\bar \gamma_1 \bar \gamma_{-1} =
\left\{
\begin{array}{l}
-q h^{-2}\\
k^{-2}
\end{array}
\right.
&\quad\!\begin{array}{l}
\mbox{for $N$ odd,} \\
\mbox{for $N$ even,}
\end{array}\\[8pt]
\bar \gamma_a \bar \gamma_{-a}=
-q k^{-2} \omega_a \omega_{a-1}
&\quad\mbox{for $a>1$,}         \nonumber
\end{array}                                             \label{bargamma}
\ee
then the elements $\bar \lambda_a$ fulfill among themselves the commutation
relations
\be
\cc{P}_{a}{}^{ab}_{cd} 
\bar \lambda_a \bar \lambda_b=0   
\ee
and the matrices $\bar e^i_a=[\bar \lambda_a,x^i]$ satisfy the
\bea
RTT-\hbox{relations:}&&  \hat R^{ij}_{kl} \bar e_a^k \bar e_b^l
= \bar e^i_c \bar e^j_d \hat R^{cd}_{ab} \label{barRee}\\
gTT-\hbox{relations:}&&
g^{ab} \bar e^i_a \bar e^j_b=g^{ij} \Lambda^{-2}, \qquad
g_{ij} \bar e^i_a \bar e^j_b=g_{ab} \Lambda^{-2},       \label{bargtt}\\
\hbox{normalization:}&& \bar e^0_0 \bar e^0_0=\Lambda^{-2}.  \label{barnor}
\eea
\end{theorem}
The conditions (\ref{ree}),(\ref{gtt}), (\ref{nor}) in theorem 2 and
(\ref{barRee}), (\ref{bargtt}), (\ref{barnor}) in theorem 4
are equivalent to the defining relations satisfied by the generators 
$\cc{L}^{\pm}{}^i_a$~\cite{FadResTak89} of $U_q^{\pm}(so(N))$, i.e.
we have found a `local realization' of the two Borel subalgebras 
$U_q^{\pm}(so(N))$ of $U_q(so(N))$. 
We can ask, under which circumstances we can `glue' them together
to construct a realization of the whole of $U_q(so(N))$.
The answer is given by the following theorem~\cite{CerFioMad00}.
\begin{theorem}                                          \label{theor5}
In the case of odd $N$ with the $\lambda_j,\bar\lambda_j$ 
defined as in (\ref{deflambda}) and (\ref{defbarlambda}) and with 
coefficients given by
\be
\begin{array}{ll}
\gamma_0 = -q^{-\frac 12}h^{-1}, & \gamma_1^2 = -q^{-2} h^{-2}, \\[4pt]
\gamma_a^2 = -q^{-2}\omega_a \omega_{a-1} k^{-2} \quad \mbox{for $a>1$,} &
\gamma_a = q \gamma_{-a}   \quad \mbox{for $a \le 1$,} \\[4pt]
\bar \gamma_a=-q \gamma_a
\end{array}                                              \label{fixgamma} 
\ee
then the matrices $e^i_a$, $\bar e^i_a$ satisfy the relations
\be
e^i_i \bar e^i_i=1 \quad \mbox{(no sum over i),} \qquad 
\hat R^{cd}_{ab}\bar e_c^i e_d^j= \hat R^{ij}_{kl}e_a^k\bar e_b^l
\label{equiv}
\ee
\end{theorem}
The $\gamma_a,\bar\gamma_a$ for $a \neq 0$ are imaginary
and fixed only up to a sign. This has as a consequence that the
homomorphism $\varphi$ does not preserve the star structure of $U_q(so(N))$.
It can be proven~\cite{CerFioMad00} that it is not possible to 
extend this theorem to the case of even $N$, because it is not possible
to verify (\ref{equiv}).

For the frames $\theta^a,\bar\theta^a$ we find:
\be
\theta^a=\theta^a_l\xi^l=
\Lambda^{-2}g^{ab}[\lambda_b,x^j]g_{jl}\xi^l, \qquad
\bar\theta^a=\bar\theta^a_l\bar\xi^l=
\Lambda^2g^{ab}[\bar\lambda_b,x^j]g_{jl}\bar\xi^l.
\ee
They commute both with the coordinates and with $\Lambda$ and the matrix 
elements $\bar\theta^a_i,\theta^a_i$ fulfill 
\bea
\hat R_{cd}^{ab} \theta^d_j \theta^c_i = \theta^b_l
\theta^a_k\hat R_{ij}^{kl}, &
\hat R_{cd}^{ab} \bar\theta^d_j \bar\theta^c_i = \bar\theta^b_l
\bar\theta^a_k\hat R_{ij}^{kl}, \\
g_{ab}\theta^b_j \theta^a_i= \Lambda^{-2} g^{ij}, \: \:
g^{ij} \theta_j^b \theta_i^a=\Lambda^{-2} g_{ab}, &
g_{ab}\bar \theta^b_j \bar \theta^a_i= \Lambda^2 g^{ij}, \: \:
g^{ij} \bar \theta_j^b \bar \theta_i^a=\Lambda^2 g_{ab}.
\eea
This implies that
\be
{\cc{P}_{s,t}}^{ab}_{cd} \theta^c \theta^d=0, \qquad
\cc{P}_{s,t}{}^{ab}_{cd} \bar\theta^c \bar\theta^d = 0.
\ee
In other words,the $\theta^a$, $\bar \theta^a$ satisfy the same 
commutation relations as the $\xi^a, \bar \xi^a$.

Finally we summarize the results found in \cite{CerFioMad00} for the metric
and the covariant derivative for each of the two calculi $\Omega(\bb{R}^N_q)$
and $\Omega^*(\bb{R}^N_q)$. In the $\theta^a$, $\bar \theta^a$ basis 
respectively the actions of $g$ and $\sigma$ are
\bea
\sigma (\theta^a\otimes \theta^b) =S^{ab}{}_{cd}\,\theta^c\otimes \theta^d,
&g(\theta^a\otimes \theta^b) = g^{ab} &\mbox{for }\Omega^*(\bb{R}^N_q),\\
\sigma (\bar\theta^a\otimes \bar\theta^b) =
\bar S^{ab}{}_{cd}\,\bar\theta^c\otimes \bar\theta^d,
&g(\bar\theta^a\otimes\bar\theta^b) = g^{ab}
&\mbox{for }\bar\Omega^*(\bb{R}^N_q).
\eea
Unfortunately, it is not possible to satisfy simultaneously the metric
compatibility condition (\ref{met-comp}) and the bilinearity condition
for the torsion (\ref{torsion}). The best we can do is
to weaken the compatibility condition to a condition of proportionality.
Then for each calculus we find the two solutions for $\sigma$:
\bea
S = q\hat R, & S = (q\hat R)^{-1} &\mbox{for }\Omega^*(\bb{R}^N_q),\\
\bar S = q\hat R, & \bar S = (q\hat R)^{-1}
&\mbox{for }\bar\Omega^*(\bb{R}^N_q).
\eea
This implies that the covariant derivatives and metric are compatible
only up to a conformal factor:
\be
S^{ae}_{df} g^{fg} S^{cb}_{eg}=q^{\pm 2} g^{ac} \delta^b_d, \qquad
\bar S {}^{ae}{}_{df} g^{fg} \bar S {}^{cb}{}_{eg} 
= q^{\pm 2} g^{ac} \delta^b_d.
\ee
In the $\xi^i$, $\bar\xi^i$ basis the actions of $g$ and $\sigma$ become
\bea
g(\xi^i \otimes \xi^j)= g^{ij} \Lambda^2 &
\sigma(\xi^i \otimes \xi^j)=S^{ij}_{hk} \xi^h \otimes \xi^k &
\mbox{for }\Omega^*(\bb{R}^N_q),\\
g(\bar \xi^i\otimes \bar\xi^j) = g^{ij} \Lambda^{-2},
&\sigma (\bar \xi^i\otimes \bar \xi^j) =
\bar S^{ij}{}_{hk}\,\bar\xi^h\otimes \bar\xi^k.
&\mbox{for }\bar\Omega^*(\bb{R}^N_q).
\eea
According to (\ref{covdev0})
the two associated covariant derivatives, one for each choice of $\sigma$, are
\be
D \xi = -\theta \otimes \xi + \sigma (\xi \otimes \theta), \qquad
\bar D \bar \xi = -\bar \theta \otimes \bar \xi + 
\sigma (\bar \xi \otimes \bar \theta).
\ee
The associated linear curvatures Curv and $\overline{\mbox{Curv}}$ vanish,
as was to be expected, because $\bb R^N_q$ should be flat.

\end{document}